\documentclass{amsart}

\usepackage{amsthm,amsmath,graphicx,epstopdf,comment,longtable,caption}

\usepackage{booktabs}
\usepackage{soul,xcolor,hyperref,supertabular,multicol,placeins}

\usepackage{floatrow,subcaption}
\usepackage[margin=4cm]{geometry}

\theoremstyle{plain}
\numberwithin{equation}{section}

\newtheorem*{problem}{Problem}

\theoremstyle{definition}

\title[Neural networks and knot invariants]{A neural network approach to predicting and computing knot invariants}

\author{Mark C. Hughes}

\begin{document}
\setstcolor{red}

\maketitle

\begin{abstract} 
In this paper we use artificial neural networks to predict and help compute the values of certain knot invariants.  In particular, we show that neural networks are able to predict when a knot is quasipositive with a high degree of accuracy.  Given a knot with unknown quasipositivity we use these predictions to identify braid representatives that are likely to be quasipositive, which we then subject to further testing to verify.  Using these techniques we identify 84 new quasipositive 11 and 12--crossing knots.  Furthermore, we show that neural networks are also able to predict and help compute the slice genus and Ozsv\'{a}th-Szab\'{o} $\tau$--invariant of knots.
\end{abstract}

\section{Introduction}

Recently developed techniques in machine learning and artificial neural networks, combined with increased computing capabilities, allow computers to recognize subtle patterns in complex data and make surprisingly accurate predictions.  Neural networks have been shown to be adept at discerning patterns in diverse data sets, and are being applied successfully to problems in image and video recognition \cite{video,image}, financial modeling \cite{finance,finance2}, natural language processing \cite{language,language2}, medical diagnostics \cite{cancer,diagnosis}, as well as many other challenging tasks.  With numerous commercial and academic applications, neural networks have been the focus of a great deal of study both in industry and academia.   

In this paper we describe an application of this technology to studying knots in the 3--sphere $S^3$.  We construct and train neural networks which are able to accurately predict the values of certain invariants that are otherwise difficult to compute.  More precisely, we construct neural networks that are able to predict whether a given knot $K$ is quasipositive, as well as the values of the slice genus $g_4(K)$ and Ozsv\'{a}th-Szab\'{o} $\tau$--invariant $\tau (K)$.  These networks predict quasipositivity with mean accuracy greater than $99.93\%$ on blind validation sets, as well as $g_4(K)$ and $\tau(K)$ with mean validation accuracies of $93.70\%$ and $99.97\%$ respectively.  Furthermore, the networks we construct learn important known relations between these invariants, and in the case of $g_4$ perform much better than baseline results obtained by combining various slice genus bounds (which achieve a mean accuracy of $71.69\%$).  

Using the predictions generated by these networks we identify braid representatives of knots that are highly favored to be quasipositive, which we then attempt to verify via a separate algorithm. Out of 179 knots with twelve crossings or less and undetermined quasipositivity, these techniques yield 84 new quasipositive knots.  We also use these techniques to compute the $\tau$--invariant for 69 new 12--crossing knots.  A similar approach may also be applied to the problems of computing $g_4(K)$ for a given knot $K$, which we illustrate below.  

In Section~\ref{sec:knotinvariants} we describe the invariants $g_4(K)$ and $\tau(K)$, as well as the notion of quasipositivity.  Section~\ref{sec:neuralnetworks} contains background material on neural networks and machine learning, while Section~\ref{sec:results} contains a description of how they are applied to model knot invariants, and is where we discuss our main results.  We include a discussion of the limitations of neural networks in this setting in Section~\ref{sec:limitations}, and an outline of further applications in Section~\ref{sec:furtherapplications}.  Finally, in Section~\ref{sec:predictions} we record the results of our main computations, as well as the slice genus and quasipositivity predictions which our networks produce for the remaining unknown small crossing knots.  Throughout this paper, all of the neural networks discussed are implemented in Python with Theano \cite{theano}, using Keras \cite{keras}.

While these initial results are very encouraging, they are only a sampling of the numerous potential applications of neural networks to knot theory.  For example, besides making predictions and aiding computations, a great deal of information could be gleaned from studying the trained networks themselves.  Any insight into how these networks make predictions may reveal valuable information about the underlying knots and their associated invariants.

\subsection*{Acknowlegdements}  The author would like to thank Eli Grigsby, Jonathan Paprocki, and David Wingate for helpful conversations.

\section{Knot invariants}
\label{sec:knotinvariants}

One useful way of studying knots is by looking at the surfaces they bound.  Suppose that $K \subset S^3$ is a knot.  We may think of $S^3$ as the boundary of $D^4$, and consider the class of orientable compact smoothly embedded surfaces $F$ in $D^4$ with $\partial F = K \subset S^3 = \partial D^4$.  Such surfaces $F$ are called \emph{slice surfaces} for the knot $K$.  Denoting the topological genus of a surface $F$ by $g(F)$, we define the (smooth) \emph{slice genus} of a knot $K$ to be 
\[
g_4(K) = \min \{ g(F) | F \text{ is a slice surface for }K\}.
\]
The problem of determining the slice genus of knots has been studied since the 1960s \cite{fox1962}.  Despite the fact that it has important implications in low-dimensional topology, $g_4(K)$ remains very difficult to determine in general.

Fortunately, the slice genus is related to another property of knots called \emph{quasipositivity}.  For each $n \geq 1$, let $B_n$ denote the $n$--strand braid group, and let $\sigma_1, \ldots , \sigma_{n-1}$ denote the standard Artin generators of $B_n$.  A \emph{band} in $B_n$ is a braid of the form $b=\alpha \cdot \sigma_j^{\pm 1} \cdot \alpha^{-1}$, where $\alpha$ is an arbitrary braid word.  Clearly any braid can be written as a product of bands.  A braid is called \emph{quasipositive} if it can be written as a product of only positive bands $b=\alpha \cdot \sigma_j \cdot \alpha^{-1}$, while a knot is called \emph{quasipositive} if it is the closure of a quasipositive braid.  

Rudolph \cite{Rudolph1993} proved that if $K$ is the closure of a quasipositive braid $\beta \in B_n$, then 
\[
g_4(K) = \frac{w(\beta) - n +1}{2},
\]    
where $w(\beta)$ is the writhe of $\beta$, i.e.\ the algebraic length of $\beta$ in the generators $\sigma_1^{\pm 1}, \ldots , \sigma_{n-1}^{\pm 1}$.  Hence if $K$ is quasipositive, the slice genus of $K$ can be readily determined provided we can find a quasipositive braid representative.

Detecting quasipositivity thus provides a means to understanding the slice genus, as well as other important features of a knot.  Unfortunately, there is currently no known algorithm for determining whether a given knot is quasipositive, nor is there a procedure for finding a quasipositive braid representative of a given quasipositive knot.  

Another related invariant of a knot $K$ is the $\tau$--invariant $\tau(K)$, defined by Ozsv\'{a}th and Szab\'{o} in \cite{ozsvathszabo}.  This integer-valued invariant is defined in terms of the Heegaard-Floer homology of $K$, and satisfies the following useful properties
\begin{enumerate}
\itemsep0.5em
\item $\tau (K_1 \# K_2) = \tau (K_1) + \tau (K_2)$,\label{tauprop1}
\item $\tau \left(\overline{K}\right) = - \tau (K)$, where $\overline{K}$ is the mirror of $K$, \label{tauprop2}
\item $| \tau (K) | \leq g_4(K)$.\label{tauprop3}
\end{enumerate}
From (\ref{tauprop1}) and (\ref{tauprop2}) it follows that $\tau (K)$ is an invariant of the smooth concordance class of $K$.  Plamenevskaya \cite{Plamenevskaya} proved that for a quasipositive knot $K$ the bound in (\ref{tauprop3}) can be strengthened to the equality $\tau (K) = g_4(K)$.   

Although Manolescu, Ozsv\'{a}th, and Sarkar proved that $\tau$ is algorithmic \cite{manolescu}, it remains computationally difficult to determine for knots with many crossings. In \cite{baldwingillam} Baldwin and Gillam compute $\tau (K)$ for all 11--crossing non-alternating knots, though for knots with 12 or more crossings the values of $\tau$ are absent from the literature.

For many knots, the challenge of computing the above invariants can often be minimized if we are able to produce accurate predictions ahead of time.  For example, the problem of verifying quasipositivity can be simplified greatly if we are able to identify, with reasonable speed and accuracy, braid representatives which are most likely to be quasipositive.  Indeed, given a knot $K$ we can generate large samples of braid representatives of $K$.  If we then have a way to select the braids which are most likely to be quasipositive, we can shrink down the search space to a few likely quasipositive candidates, which could be checked by hand or computer algorithm.  Similarly, when trying to determine the slice genus of a knot $K$ it is often helpful to have an estimate for the value of $g_4(K)$ as a starting point.  An accurate estimate of $g_4(K)$ often indicates which approach is most likely to be successful and should be tried first, i.e.\ whether to try finding obstructions or look for surfaces which explicitly realize the estimated value.  In Section~\ref{sec:neuralnetworks} we construct neural networks which provide these predictions.  Of course, such predictions will always need to be verified by a rigorous proof, but for practical purposes they are useful as they are often able to tell us what we should be trying to prove.

\section{Neural networks}
\label{sec:neuralnetworks}

The topic of neural networks has been the focus of a great deal of study and interest in the past 25 years.  While they are useful in a variety of contexts for numerous applications, we will use them to study the following interpolation-type problem:

\begin{problem}
Let $\mathcal{K}$ be a set, with $f:\mathcal{K} \rightarrow \mathbb{R}^m$ a function.  Suppose that the values of $f$ are known only on a proper finite subset $\mathcal{S} \subset \mathcal{K}$.  Can we predict the values of $f(x)$ for $x \in \mathcal{K} \backslash \mathcal{S}$?
\end{problem}

Naturally if we want any hope of succeeding at this problem, the set $\mathcal{K}$ will need to have some additional structure which is related to the function $f$.  By observing the values of $f$ on $\mathcal{S}$, a neural network seeks to detect this relevant structure and deduce its relation to the function $f$.

It is useful to note that the function $f$ could be used to describe a number of different features of the set $\mathcal{K}$.  For example, suppose that each of the elements of $\mathcal{K}$ can be sorted into one of $k$ distinct classes $\mathcal{C}_1, \ldots, \mathcal{C}_k$. Such a classification can be represented as a function $f:\mathcal{K} \rightarrow \mathbb{R}^k$, where $f$ sends elements of $\mathcal{C}_j$ to the $j$th standard basis vector of $\mathbb{R}^k$.  Interpolating the function $f$ on the set $\mathcal{K}\backslash \mathcal{S}$ then corresponds to predicting the classes of elements with unknown classification.  These problems are thus aptly named \emph{classification problems}.  On the other hand, often the function $f$ is continuous and we would like to predict it's values on $\mathcal{K} \backslash \mathcal{S}$.  Such problems are called \emph{regression problems}.

Each of these types of problems can be approached via neural networks, which we describe below.  Although the term \emph{neural network} refers to a wide class of network type functions with varying sophistication, we only describe the simplest form here, called a \emph{feed-forward neural network}.

Let $n_0, \ldots , n_q \in \mathbb{N}$.  For each integer $1 \leq j \leq q$, let $h_j:\mathbb{R} \rightarrow \mathbb{R}$ be a function.  Furthermore, let $A_j$ be an $n_j \times n_{j-1}$ real-valued matrix, let $b_j \in \mathbb{R}^{n_j}$, and let $T_j: \mathbb{R}^{n_{j-1}} \rightarrow \mathbb{R}^{n_j}$ be the affine function
\[
T_j(x) = A_j x + b_j.
\]
Then the function $\psi:\mathbb{R}^{n_0} \rightarrow \mathbb{R}^{n_q}$ defined by the composition
\[
\mathbb{R}^{n_0} \xrightarrow{\text{ } T_1 \text{ }} \mathbb{R}^{n_1} \xrightarrow{\text{ } h_1 \text{ }} \mathbb{R}^{n_1} \xrightarrow{\text{ } T_2 \text{ }} \mathbb{R}^{n_2} \xrightarrow{\text{ } h_2 \text{ }} \mathbb{R}^{n_2}\xrightarrow{\text{ } T_3 \text{ }} \cdots \xrightarrow{\text{ } T_q \text{ }}\mathbb{R}^{n_q} \xrightarrow{\text{ } h_q \text{ }} \mathbb{R}^{n_q}
\]
is called a \emph{feed-forward neural network}.  Here, the function $h_j$ are called \emph{activation functions}, and act on $\mathbb{R}^{n_j}$ component-wise.  Obviously the network $\psi$ depends on the values of the weight matrices $A_j$ and bias vectors $b_j$, though we suppress these dependencies in the notation for convenience. 

It is often convenient to think of the network $\psi$ as being divided into a sequence of \emph{layers}.  The input vectors in $\mathbb{R}^{n_0}$ form the \emph{input layer}.  For each $1 \leq j \leq q-1$ the output of $h_j$ forms the $j$th \emph{hidden layer}, while the output of $h_q$ is called the \emph{output layer}.  The name neural network is due to the network structure that can be seen by considering how the components of each layer are obtained from the previous layer.  Indeed, the components of the hidden and output layers are each obtained as a weighted sum of the components of the previous layer and bias term $b_j$, with activation function $h_j$ applied.  Each component can be thought of as a single neuron in the network, taking input from each of the neurons in the previous layer, and outputting a signal which has been modulated by its activation function.    

In practice the activation functions $h_j$ on the hidden layers are typically chosen from the set
\[
\left\{ \tanh(x), \max \{0,x\}, \left(1+e^{-x}\right)^{-1} \right\}.
\]
They can be thought of as imposing a threshold which decides when the corresponding neuron `fires' depending on the input value $x$.  For regression problems, we often take the final activation function $h_q$ to be the identity.  For classification functions, instead of applying an activation function $h_q$ component-wise to the output of $T_q$, we apply the \emph{softmax} function 
\[
\sigma (x_1, \ldots , x_{n_q}) = \left( \frac{\exp(x_1)}{\sum_{j=1}^{n_q} \exp (x_j)}, \ldots , \frac{\exp(x_{n_q})}{\sum_{j=1}^{n_q} \exp (x_j)}\right).
\]
In this case the output vector $(p_1, \ldots ,p_{n_q})$ satisfies $\sum_{j=1}^{n_q} p_j = 1$, where $p_j$ can be interpreted as the probability that the input vector represents an element living in the $j$th class $\mathcal{C}_j$.

While the function type described above may seem restrictive, neural networks actually suffice to approximate arbitrary continuous functions on compact subsets of $\mathbb{R}^{n_0}$. Indeed, it is known that neural networks with as few as one hidden layer and arbitrary nonconstant bounded activation functions can be used to approximate any continuous function on a compact subset of $\mathbb{R}^{n_0}$ (see, e.g. \cite{Hornik1991}).

When using a neural network to interpolate a function $f:\mathbb{R}^{n} \rightarrow \mathbb{R}^{m}$, the network architecture is selected first. This includes selecting the number and sizes of the hidden layers, as well as the activation functions  $h_j$.  The weight matrices $A_j$ and bias vectors $b_j$ are randomly initialized, and tuned via a process called \emph{training}.  

To train the network, we first select an error function.  For regression problems, this error function is often chosen to be the \emph{mean-squared error}
\[
E_{MS}(A_1, \ldots, A_q, b_1 , \ldots , b_q) = \sum_{x\in \mathcal{S}} \| f(x)-\psi(x) \|^2.
\]
Here, we think of $E_{MS}$ as being a function of the weight matrices $A_j$ and bias vectors $b_j$.  For classification problems, we often use the \emph{cross-entropy error} function, defined as
\[
E_{CE}(A_1, \ldots, A_q, b_1 , \ldots , b_q)=-\sum_{x \in \mathcal{S}}\sum_{j=1}^{n_q}f_j(x)\ln(\psi_j(s)),
\]
where $f=(f_1,\ldots ,f_{n_q})$ and $\psi=(\psi_1,\ldots ,\psi_{n_q})$.

Training the network $\psi$ then proceeds by some variant of gradient descent for the error function $E$.  In the simplest approach, the gradient $\nabla E$ is computed with respect to the weight matrices $A_j$ and bias vectors $b_j$, and the values of these parameters are adjusted in the direction of $-\nabla E$.  This process is repeated a number of times in an attempt to minimize $E$.  Each time, the error function $E$ and its gradient can be computed on either individual points in $\mathcal{S}$, or on batches of points.

Throughout the procedure of constructing and training of the neural network $\psi$, there are a number of choices that must be made.  From network architecture parameters (including the number and dimensions of the hidden layers), to the choice of error function, to the training algorithm and learning parameters used, each choice affects how effective the network is at learning and approximating the function $f$.  While there are general guidelines which suggest appropriate choices for these settings, much of the work in constructing an accurate predictor lies in testing different choices to see which yield the most accurate predictions.  To accomplish this, a set of hold-out data, called a \emph{test set}, is removed from the set $\mathcal{S}$. The test data is not involved in the training of $\psi$, but is used to test the accuracy of the resulting network $\psi$ after training.  

\section{Predicting knot invariants}
\label{sec:results}

We now describe how these techniques are applied to the problem of predicting and computing knot invariants.  We focus on the set of knots with 12 or fewer crossings, of which there are 2977, using knot data obtained from KnotInfo \cite{knotinfo}. 

\subsection*{Data preparation} To begin, we must select a way to identify these knots with a subset of some $\mathbb{R}^n$.  Each of these knot can be represented as a braid word of length 19 or less in the Artin generators $\sigma^{\pm 1}_1, \ldots, \sigma^{\pm 1}_6$,  which we identify with a $228$--tuple of numbers using a process called \emph{one-hot encoding}.  To do this we think of each braid as being divided into 19 slots, each of which can either be empty (corresponding to no crossing), or hold a single letter $\sigma_j^{\pm 1}$.  To each of these 19 slots we assign a separate copy of $\mathbb{R}^{12}$, denoted $\mathbb{R}^{12}_j$ for $1 \leq j \leq 19$, and we identify the contents of the $j$th slot with a vector in $\mathbb{R}_j^{12}$.  Empty slots are identified with the origin, while slots containing one of the twelve letters $\sigma_1, \ldots, \sigma_6, \sigma_1^{-1}, \ldots , \sigma_6^{-1}$ are identified with one of the twelve standard basis vectors in $\mathbb{R}^{12}_j$ respectively.  More precisely, if the $j$th slot contains the letter $\sigma_k^{\varepsilon}$, then it is identified with the $(k+3(1-\varepsilon))$th basis vector in $\mathbb{R}^{12}_j$.  This gives us a way to identify each of the 19 possible crossings with a vector in $\mathbb{R}^{12}$, and hence we can identify any braid word of length 19 or less with a vector in $\bigoplus_{j=1}^{19} \mathbb{R}_j^{12} \cong \mathbb{R}^{228}$.  While there are certainly more compact ways to encode braid words as vectors in some Euclidean space, the procedure described here ensures that we are not introducing any unintentional ordering relations on the braids we represent.

By increasing our representation space from $\mathbb{R}^{228}$ to $\mathbb{R}^{240}\cong \mathbb{R}^{228} \oplus \mathbb{R}^{12}$ we also encode twelve additional features which describe invariants or properties of the associate knots.  Five binary variables are introduced to describe whether the corresponding knot is alternating, fibered, a positive braid closure, or large or small, together with integer-valued variables which encode the crossing number, Seifert genus, braid index, signature, arc index, determinant, and Rasmussen invariant of the knot (see KnotInfo for descriptions of each of these invariants).  Each variable is then normalized so that it has zero mean and unit standard deviation.  While this normalization is not necessary in theory, in practice it greatly decreases the time needed to train the network effectively.

We included the above features in our data set not because of any known or suspected relations to quasipositivity, the slice genus, or $\tau$, but rather because their values are completely listed in KnotInfo for all of the knots we study.  Indeed, an interesting direction of study would be to test what effect the inclusion of different variables has on the ability of a neural network to predict seemingly unrelated invariants.  Our choice to include only invariants which have been completely computed is not strictly necessary, as there are various techniques for imputing missing data in the machine learning literature (e.g., \cite{imputing}). 

Thus far we have fixed a single braid word representative for each knot we are considering, which is overly restrictive.  We want the networks to detect topological invariants of knots, instead of simply memorizing the features of a select set of braid words.  By modifying the braid representatives in our collection by random Markov and Reidemeister moves, we generate additional braid representatives for each knot and append them to our data set.  In the end we generate 32 different braid representatives for each knot, yielding a total of 95,264 data vectors that can be used to train our networks.    

One final step of data preparation is needed before we proceed to define and train our neural networks.  For each of the features of interest (quasipositivity, slice genus, or the $\tau$--invariant), the data is split into three sets. The \emph{training data} consists of the majority of the data vectors for which the value of the specified feature is known.  It is the data that is used during the training procedure to determine appropriate values for the weight matrices $A_j$ and bias vectors $b_j$.  From this data we separate a smaller set of \emph{test data}, for which the values of the specified feature are known, but which are not used during training.  This is the data we use to test the accuracy of the trained network.  The remaining data corresponds to knots for which we do not know the value of the desired feature, but would like to predict via the trained neural network.

Because we want to ensure that the neural networks are actually predicting knot invariants and not just recognizing Reidemeister and Markov moves, when we include one braid representative of a given knot $K$ in the test data set we make sure to include all other representatives of $K$ in the test data as well.  Thus, the networks are only tested on genuinely new knot classes which they have not been exposed to during the training phase.

\subsection*{Quasipositivity detection}  We begin by describing a neural network $\psi_Q$ which attempts to detect quasipositivity in knots.  This is a classification problem, where each knot lies in one of two classes: quasipositive or not quasipositive.  As is the case with all of the networks described here, it will have an input layer of 240 neurons due to the method we used to encode the data as vectors in $\mathbb{R}^{240}$.  The input layer is followed by two hidden layers, with 500 neurons each and activation functions $h(x) = \max \{ x, 0\}$.  The output layer consists of a single neuron with sigmoidal activation function $\sigma(x) = (1+e^{-x})^{-1}$ (which is a specialization of the softmax function to classification problems with only two possible classes).  

At each hidden layer we also apply a procedure known as \emph{Dropout} \cite{dropout}, which is a process by which certain neurons are ignored randomly during training, and the resulting thinned networks are averaged.  This prevents the neural network from being ``over-fit" to the training data, a problem which limits the ability of the network to make predictions on new data.

The network $\psi_Q$ is trained using a variant of gradient descent called \emph{Adam} \cite{adam}, with the cross-entropy function as our error function.  The entire training data set is fed through the algorithm a total of 26 times (in other words, 26 \emph{epochs}) to yield a trained network $\psi_Q$.  This entire procedure is repeated five times, creating five different neural networks $\psi_Q^1, \ldots, \psi_Q^5$.  Each time different training and test data sets are used, where the test data sets consist of 1600 braid representatives each, or 500 distinct knot types.

Given an input vector $v \in \mathbb{R}^{240}$ representing some knot, the output of each $\psi^j_Q$ is a value $\psi^j_Q(v) \in [0,1]$, which can be interpreted as the network-assigned probability that the vector $v$ corresponds to a quasipositive knot.  To make concrete predictions, we will set a cutoff value and predict that knots with $\psi^j_Q(v) \geq 0.5$ are quasipositive, while knots with $\psi^j_Q(v) < 0.5$ are not quasipositive.  On average the networks $\psi^j_Q$ are able to correctly predict quasipositivity on the test data sets $99.93\%$ of the time.   

Once we have estimated the accuracy of the networks $\psi^j_Q$, we train one final network $\psi_Q$, this time using all available data for training.  
Applying $\psi_Q$ to the remaining data vectors with unknown quasipositivity allows us to quickly identify braid representatives which are likely to be quasipositive.  These candidate braids are then subjected to a separate algorithm which searches for positive band decompositions.  

Using these techniques we find 72 new quasipositive knots with 11 and 12 crossings.  Furthermore, using similar methods we find an additional 12 knots which are \emph{quasinegative}, i.e.\ knots whose mirror images are quasipositive.  These results are summarized in Section~\ref{sec:predictions}.  Because the $\tau$--invariant and slice genus of quaspositive knots are easily computed given a quasipositive braid representative (see Section~\ref{sec:knotinvariants}), we also include the values of $\tau(K)$ and $g_4(K)$ for each of the knot types listed.  The values of $\tau (K)$ for these knots are not currently listed in KnotInfo, though Baldwin and Gillam \cite{baldwingillam} do compute $\tau (K)$ for all 11--crossing knots (using differing sign conventions).   A benefit to the quasipositive braid presentations in Section~\ref{sec:predictions} is that they allow us to explicitly construct slice surfaces $F \subset D^4$ with $\partial F = K$, which realize $g_4(K)$ (see \cite{Rudolph1983}).

\subsection*{Slice genus and Ozsv\'{a}th-Szab\'{o} $\tau$--invariant detection}  We now describe neural networks $\psi_g$ and $\psi_\tau$ which predict the slice genus and $\tau$--invariant respectively.  Because both of these invariants take numeric values, we treat the modeling of them as regression problems.  Note however, that since only a finite range of values for $g_4$ and $\tau$ show up among smaller crossing knots, we could approach both as finite classification problems, though we find that the regression approach yields marginally better results.

Each of the networks $\psi_g$ and $\psi_\tau$ have two hidden layers, with activation functions $h(x)=\max \{0,x\}$, and single neuron output layer with identity activation function.  Each of the hidden layers in $\psi_g$ has 500 neurons, while the hidden layers of $\psi_\tau$ have 750 neurons each.  We employ Dropout to both networks to avoid over-fitting, and train them both using the Adam algorithm with mean-squared error loss function.  The network $\psi_g$ is trained for 27 epochs (i.e.\ the training data set is passed through the training algorithm 27 times), while $\psi_\tau$ is trained for 37 epochs.  As was the case above with detecting quasipositivity, the networks $\psi_g$ and $\psi_\tau$ are initialized and trained five separate times on different training data sets, and achieve mean accuracies on the corresponding test data sets of $93.70\%$ and $99.97\%$ respectively.

During the preparation of this paper, the information in KnotInfo was updated to reflect McCoy's results in \cite{Mccoy}, where he computes the slice genus of 615 new knots with 11 and 12 crossings. McCoy uses a computer to search for genus one cobordisms between knots, yielding bounds that suffice to compute $g_4(K)$ for many knots $K$.  These new computations provide an additional set of test data for us to compare our predictions to.  Perhaps surprisingly, on these 615 knots the network $\psi_g$ performs better than on the earlier test data sets, predicting $g_4(K)$ correctly for 596 of the 615 knots ($96.91\%$ accuracy).

It is interesting to see how well the neural network models known relations between the various invariants.  For example, consider the relation $g_4(K) \geq |\tau (K)|$.  If $x \in \mathbb{R}$, let $[x] \in \mathbb{Z}$ be the nearest integer to $x$ (rounding up when $x$ is a multiple of 0.5).  Then for each of our 95,264 data vectors $v \in \mathbb{R}^{240}$ in our data set, we have
\begin{equation}
\label{eqn:relation}
[\psi_g(v)] \geq | [ \psi_\tau (v) ]|.
\end{equation}
The relation in (\ref{eqn:relation}) remains true for all $v$ when we replace either $[\psi_g(v)]$ or $[\psi_\tau (v)]$ with a known value of $g_4(K)$ or $\tau(K)$ respectively, where $K$ is the knot type corresponding to the vector $v$.  We also have that 
\[
[\psi_g(v)] \geq \frac{|  s (K) |}{2}
\]
for all $v$, where $s$ is the Rasmussen $s$--invariant from Khovanov homology \cite{rasmussen} which satisfies $g_4(K) \geq \frac{1}{2}|s(K)|$ for all knots $K$.


\section{Limitations}
\label{sec:limitations}

An important part of designing and applying neural networks to knot theory is to understand their limitations.  Perhaps the most glaring such limitation is that although neural networks can provide helpful guidance and predictions when approaching a problem, these predictions do not constitute a proof, regardless of how confident the network seems.  Furthermore, neural networks are essentially ``black-boxes," in that is difficult to understand why they make the predictions they do.  

Another potential issue is the possibility of disparities between the data the networks train on and the data we would like to make predictions on.  Indeed, the networks are being trained on knots for which the values of a certain invariant are known, but making predictions on knots for which the values of this invariant are unknown.  It is possible that there are fundamental differences which separate this latter class of knots from the former, and which may even contribute to the apparent difficulty in computing the given invariant on these knots.  It is naive to expect that all knots will satisfy the patterns we observe among their simplest representatives.  Furthermore, even when neural networks do successfully identify patterns in the data which generalize outside of the training set, it is not guaranteed that they have learned information that we are actually interested in.  As such, neural networks should be viewed as guides that are privy to observations we may not be, but which are in no wise infallible.  

Finally, neural networks do not seem to be equally well-suited to modeling all invariants.  For example, Jonathan Paprocki has observed that neural networks (with architecture similar to the ones used above) seem to have difficulty predicting the Arf invariant of knots.  This is perhaps surprising given the fact that the Arf invariant is not overly difficult to compute for a given knot.  It would be interesting to understand what factors separate invariants that can be successfully modeled by neural networks of a given size from those which cannot.

\section{Additional applications}
\label{sec:furtherapplications}

While we have applied neural networks to find quasipositive braid representatives of knots, this is only a small example of the many potential uses for neural networks in knot theory.  Here we discuss others, focusing on the problem of computing the slice genus.   

As mentioned in Section~\ref{sec:results}, during preparation of this paper the knot tables in KnotInfo were updated to reflect results in \cite{Mccoy}.  Prior to this we had computed the slice genus of 24 new knots, all of which were included in \cite{Mccoy} but were unknown to the author at the time.  Of these 24 knots 16 are quasipositive, and thus their slice genus can be deduced from their quasipositive braid representatives in Section~\ref{sec:predictions} (from which explicit minimal genus slice surfaces can also be constructed).

Besides computing the slice genus of quasipositive knots, we can also use the network $\psi_g$ directly to make guiding predictions for other slice genus computations.  For example, in Table~\ref{tab:nonslice} we present a collection of non-slice knots, i.e.\ knots which are known to not bound slice disks.  For each knot $K$ in Table~\ref{tab:nonslice} we also include the mean of $\psi_g(v)$ for all data vectors $v$ representing the knot type $K$.  We let $\psi_g(K)$ denote this mean prediction, and order the knots in Table~\ref{tab:nonslice} by the difference $|\psi_g(K)-1|$.  Since we know that $g_4(K) \geq 1$ for each knot $K$ in Table~\ref{tab:nonslice}, we can attempt to prove that $g_4(K)=1$ by looking for genus one slice surfaces $F$ with $\partial F = K$.  As the knots near the top of Table~\ref{tab:nonslice} are strongly favored to have $g_4(K)=1$, they provide a promising list of examples to begin such a search with.

\begin{table}
\begin{tabular}{  c  c  c   }
\toprule
{\sc Knot type} & {\sc Mean $g_4(K)$ prediction ($\psi_g(K)$)} & \hspace{0.5cm} $| \psi_g(K) - 1|$ \hspace{0.5cm} \\ 
\toprule
$12n_{307}$ & $1.000014$ & $0.000014$ \\ \hline
$11n_{119}$ & $1.000032$ & $0.000032$ \\ \hline
$12a_{769}$ & $0.999919$ & $0.000081$ \\ \hline
$11n_{115}$ & $1.000197$ & $0.000197$ \\ \hline
$11a_{297}$ & $1.000241$ & $0.000241$ \\ \hline
$11a_{315}$ & $1.000332$ & $0.000332$ \\ \hline
$12a_{706}$ & $1.000378$ & $0.000378$ \\ \hline
$12a_{89}$ & $1.000387$ & $0.000387$ \\ \hline
$11a_{251}$ & $1.000395$ & $0.000395$ \\ \hline
$11a_{119}$ & $0.999599$ & $0.000401$ \\ \hline
$12a_{668}$ & $1.000420$ & $0.000420$ \\ \hline
$12n_{805}$ & $1.000518$ & $0.000518$ \\ \hline
$11a_{37}$ & $1.000562$ & $0.000562$ \\ \hline
$\vdots$ & $\vdots$ & $\vdots$ \\ \bottomrule
\end{tabular}
\caption{Non-slice knots and their slice genus predictions}
\label{tab:nonslice}
\end{table}

Proceeding with this approach, we find genus one slice surfaces for 8 of the top 13 knots in Table~\ref{tab:nonslice} (see Section~\ref{sec:predictions} for diagrams of these slice surfaces).  While the slice genus of these knots are determined independently in \cite{Mccoy}, a similar approach could be used to identify promising candidates for slice genus computations among higher crossing knots.

While these techniques can be used to identify knots that are likely amenable to computations, more careful thought is needed if we wish to target specific examples.  This could be accomplished, for example, by using $\psi_g$ to test for potential concordance relations among knots. More precisely, let $K\#J$ denote the connected sum of knots $K$ and $J$, and let $-J$ denote the mirror image of $J$ with reverse orientation.  We say that two knots $K$ and $J$ are \emph{concordant} if $K \# -J$ bounds a slice disk in $D^4$.  The set of all concordance classes forms a group (with operation the connected sum) which has been studied since it was introduced by Fox and Milnor \cite{fox1966} in 1966.  Despite its long history, however, the knot concordance group is still far from well-understood.  

It is straight-forward to show that if $K$ and $J$ are concordant, then $g_4(K) = g_4(J)$.  Thus given a knot $K$ with unknown slice genus, we could use a neural network to identify knots $J$ where $K\#-J$ is likely to be slice, and where $g_4(J)$ is known.  Finding an explicit slice disk for $K\#-J$ would then verify that $g_4(K)=g_4(J)$.  To be effective in practice, we would likely need to introduce more sophisticated network architecture and increase the size of our set of training data.  Indeed, by randomly generating other representatives from the concordance classes of $K$ and $J$ respectively, we could test multiple pairs of knots from each class to increase our chance of finding new concordance relations. As several other important knot invariants are also concordance invariants, this approach could be targeted to study the invariants of particularly difficult examples.  These techniques will be developed further by the author in a future paper.

Beyond using neural networks to compute invariants and find concordance relations, a great deal of information could also be obtained from studying the trained networks themselves.  For example, interesting questions could be asked about the networks' decision boundaries (which could be studied via their persistent homology) or the spectra of their weight matrices.  Understanding how these networks make predictions may reveal valuable new information about the underlying knots and their invariants.

\section{Predictions and Results}
\label{sec:predictions}

In this final section we collect results mentioned in previous sections, and present predictions for remaining open cases.  We begin with Table~\ref{tab:quasipositive}, which contains quasipositive braid representatives for 72 new knots.  For notational convenience, we represent them as braid words in the letters $A, B, C, D, a, b, c, d$, which represent $\sigma_1,\sigma_2,\sigma_3,\sigma_4,\sigma^{-1}_1,\sigma^{-1}_2,\sigma^{-1}_3$, and $\sigma^{-1}_4$ respectively.  The column {\sc Band centers} describes a positive band decomposition of the given braid, by specifying the locations of the center of each positive band (recall that a positive band is a braid word of the form $\alpha \sigma_j \alpha^{-1}$ for some Artin generator $\sigma_j$ and arbitrary braid word $\alpha$).  From these locations we can easily construct a positive band decomposition.  For example, the braid word $A A b A b A C B B C C$ has a positive band decomposition with bands centered in positions 1, 2, 4, 6, 7, 10, and 11.  From this we can construct the positive band decomposition
\[
(A) (A) (b A B) (b b A B B) (b b C B B) (C) (C).
\]
Table~\ref{tab:quasinegative} is similar, though it contains quasinegative knots with quasinegative braid representatives.

\begin{figure*}[h!]
    \centering
    \begin{subfigure}[t]{0.5\textwidth}
        \centering
        \includegraphics[height=0.53\textwidth]{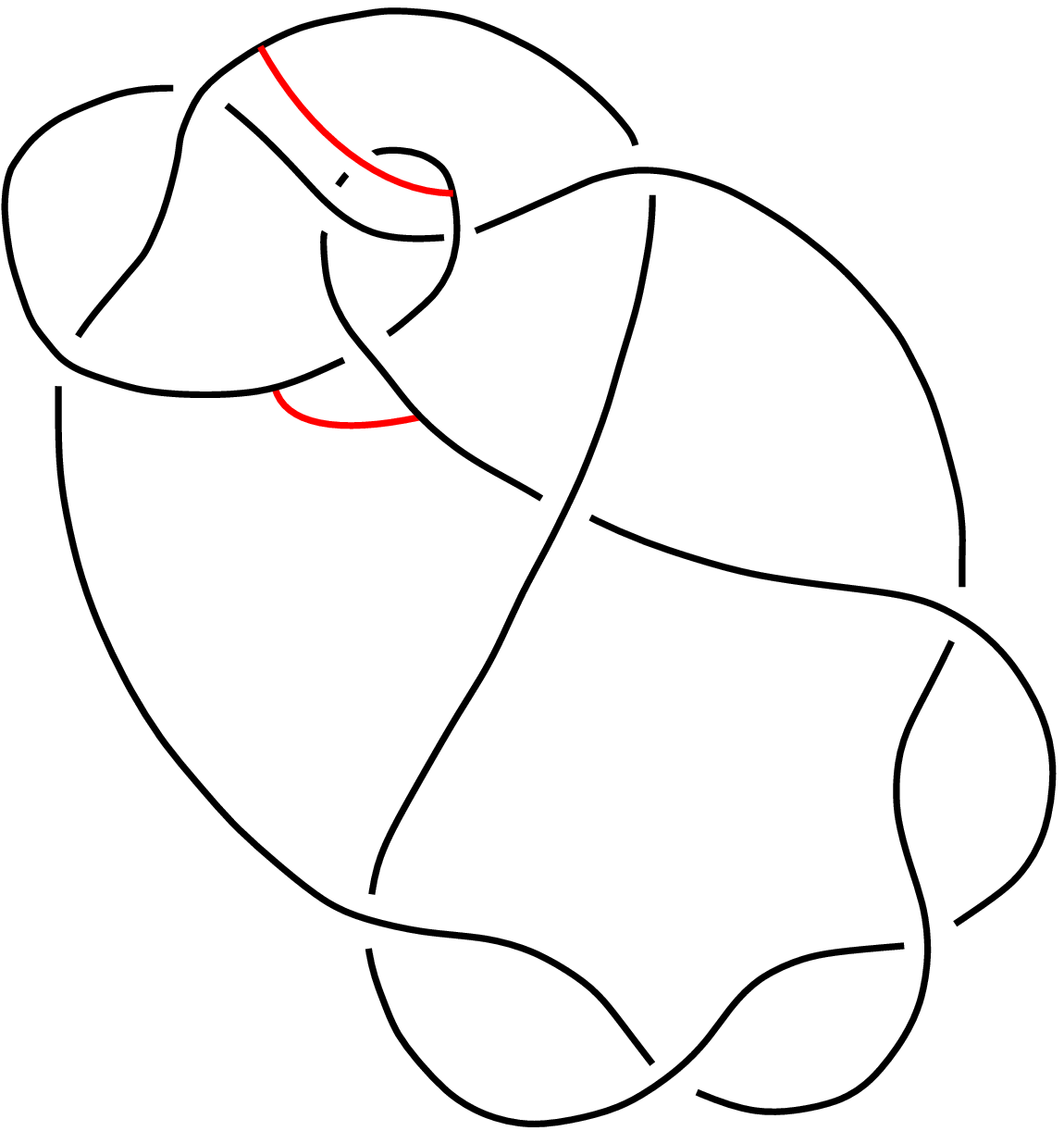}
        \caption{$11a_{37}$}
    \end{subfigure}%
    ~ 
    \begin{subfigure}[t]{0.5\textwidth}
        \centering
        \includegraphics[height=0.53\textwidth]{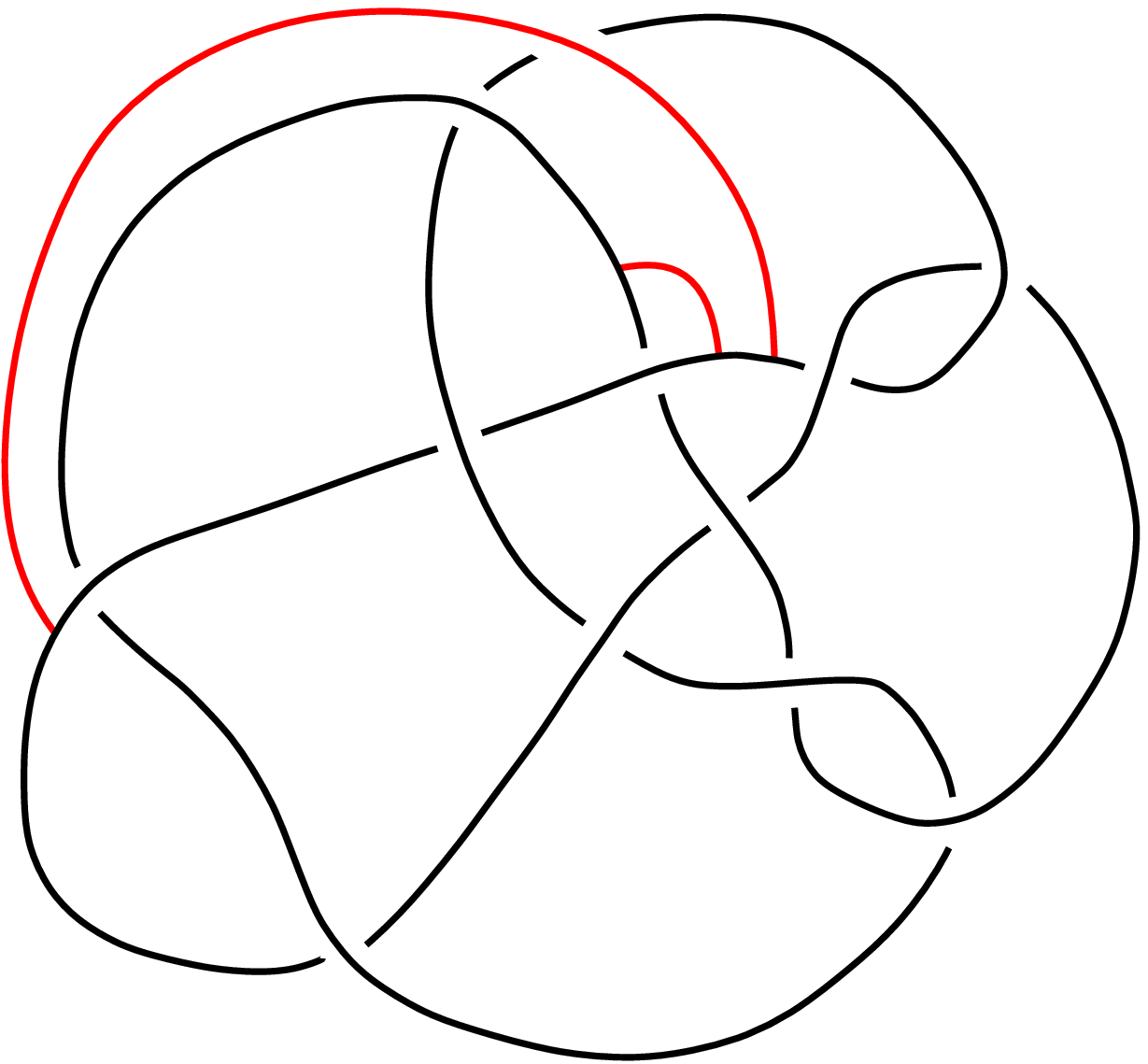}
        \caption{$11a_{315}$}
    \end{subfigure}
    
    \vspace{1cm}
        \begin{subfigure}[t]{0.5\textwidth}
        \centering
        \includegraphics[height=0.53\textwidth]{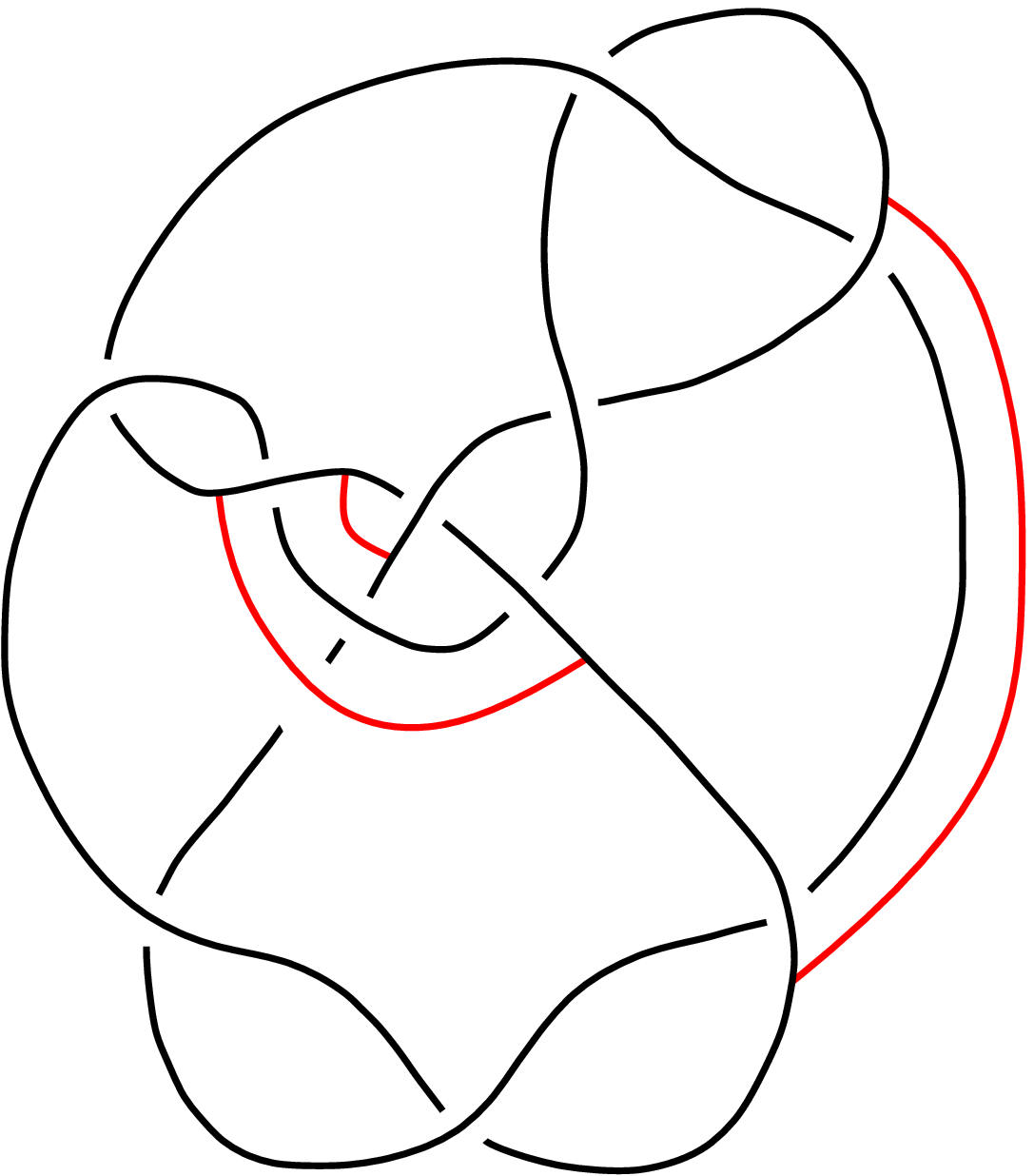}
        \caption{$11n_{115}$}
    \end{subfigure}%
    ~ 
    \begin{subfigure}[t]{0.5\textwidth}
        \centering
        \includegraphics[height=0.53\textwidth]{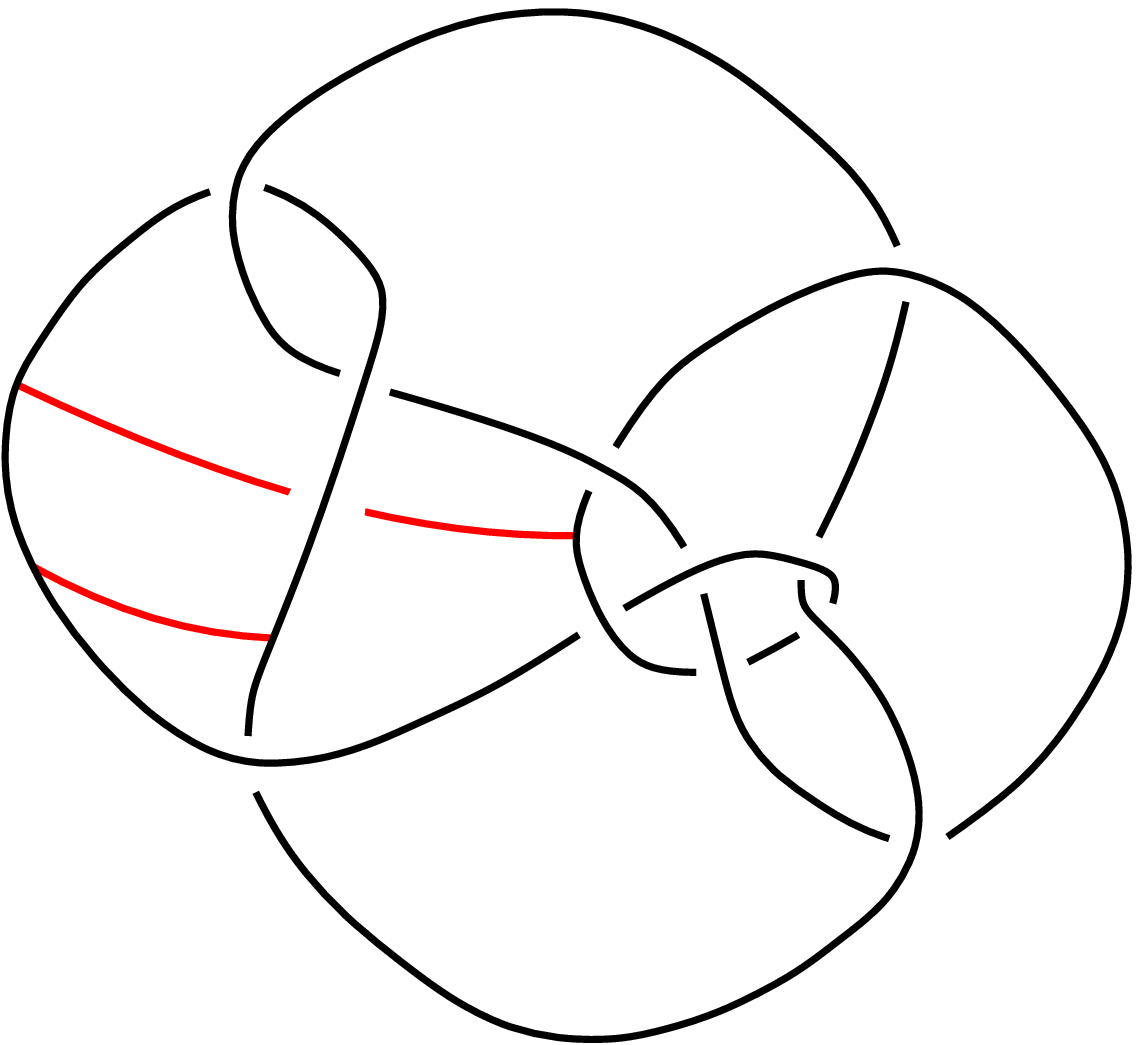}
        \caption{$11n_{119}$}
    \end{subfigure}
    
    \vspace{1cm}
        \begin{subfigure}[t]{0.5\textwidth}
        \centering
        \includegraphics[height=0.53\textwidth]{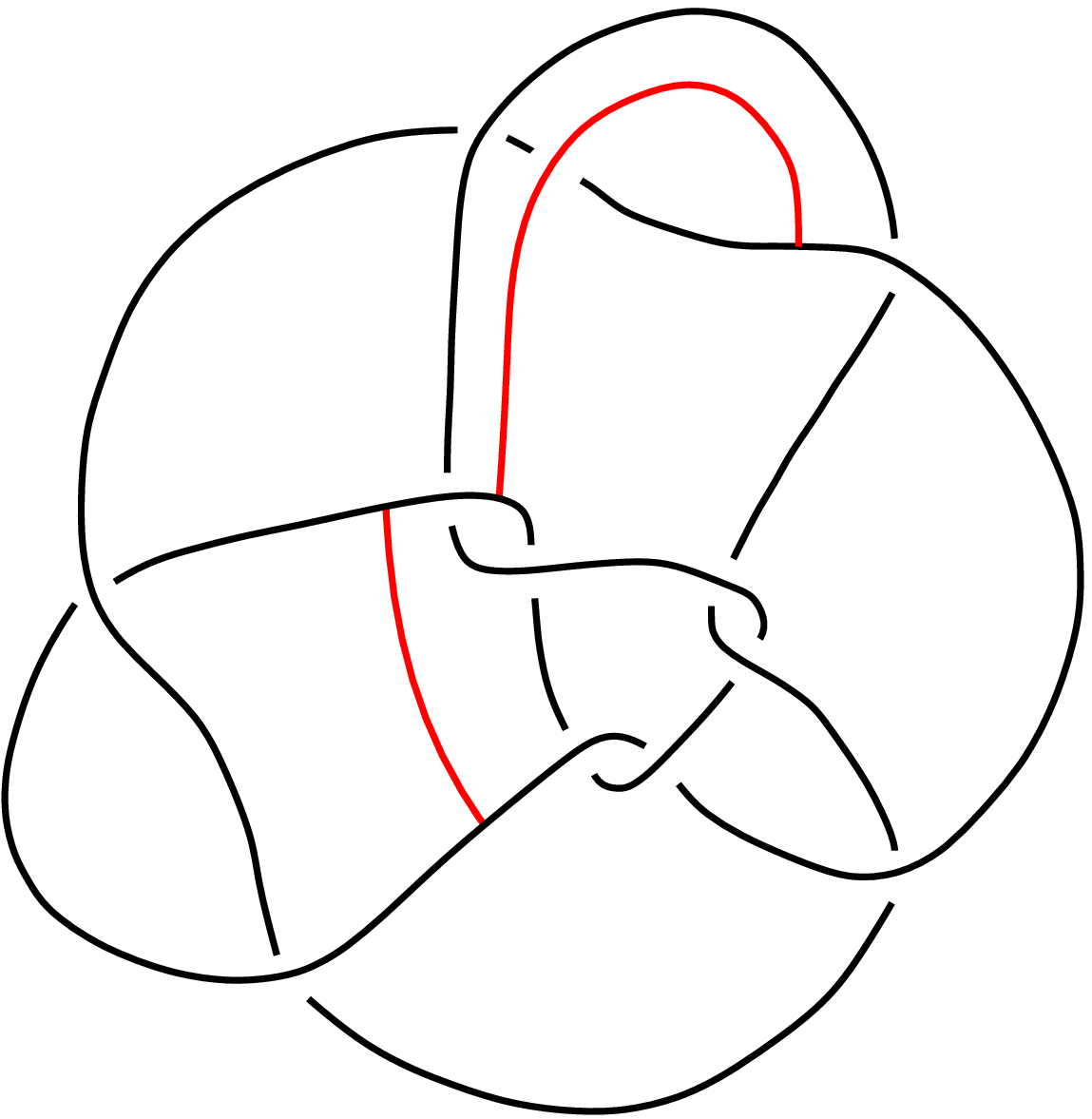}
        \caption{$11n_{179}$}
    \end{subfigure}%
    ~ 
    \begin{subfigure}[t]{0.5\textwidth}
        \centering
        \includegraphics[height=0.53\textwidth]{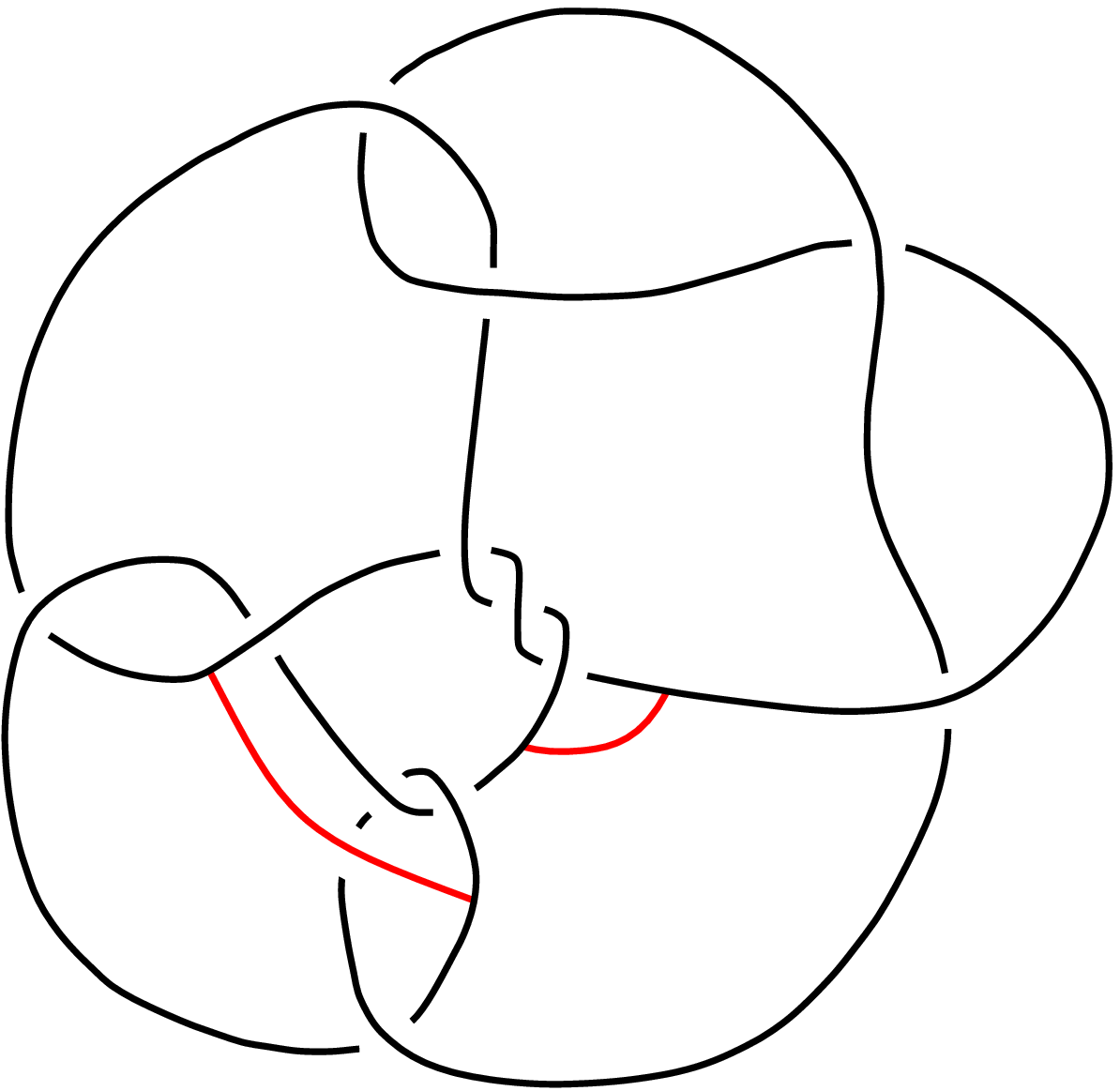}
        \caption{$12a_{89}$}
    \end{subfigure}
    
        \vspace{1cm}
        \begin{subfigure}[t]{0.5\textwidth}
        \centering
        \includegraphics[height=0.53\textwidth]{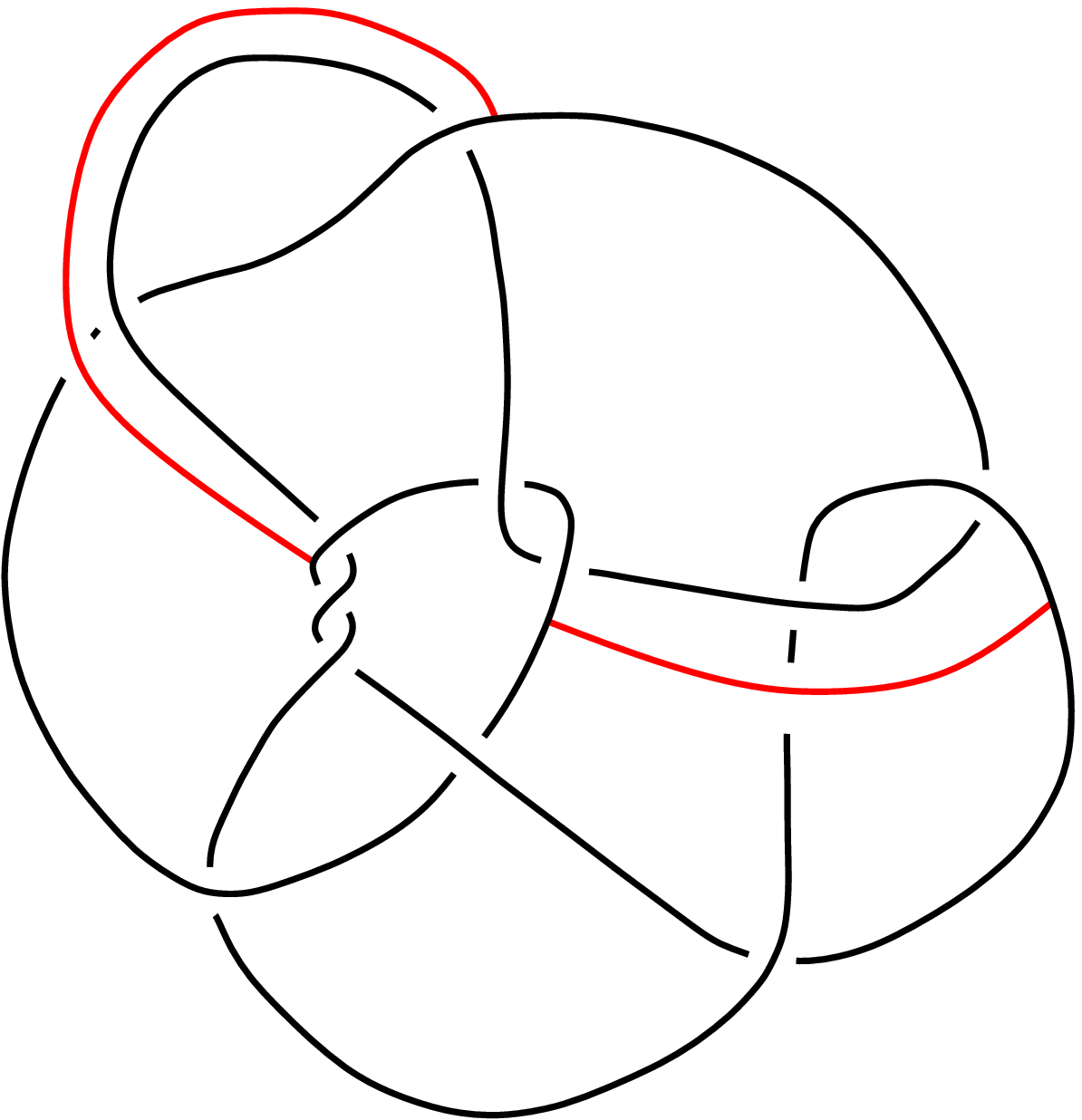}
        \caption{$12a_{769}$}
    \end{subfigure}%
    ~ 
    \begin{subfigure}[t]{0.5\textwidth}
        \centering
        \includegraphics[height=0.53\textwidth]{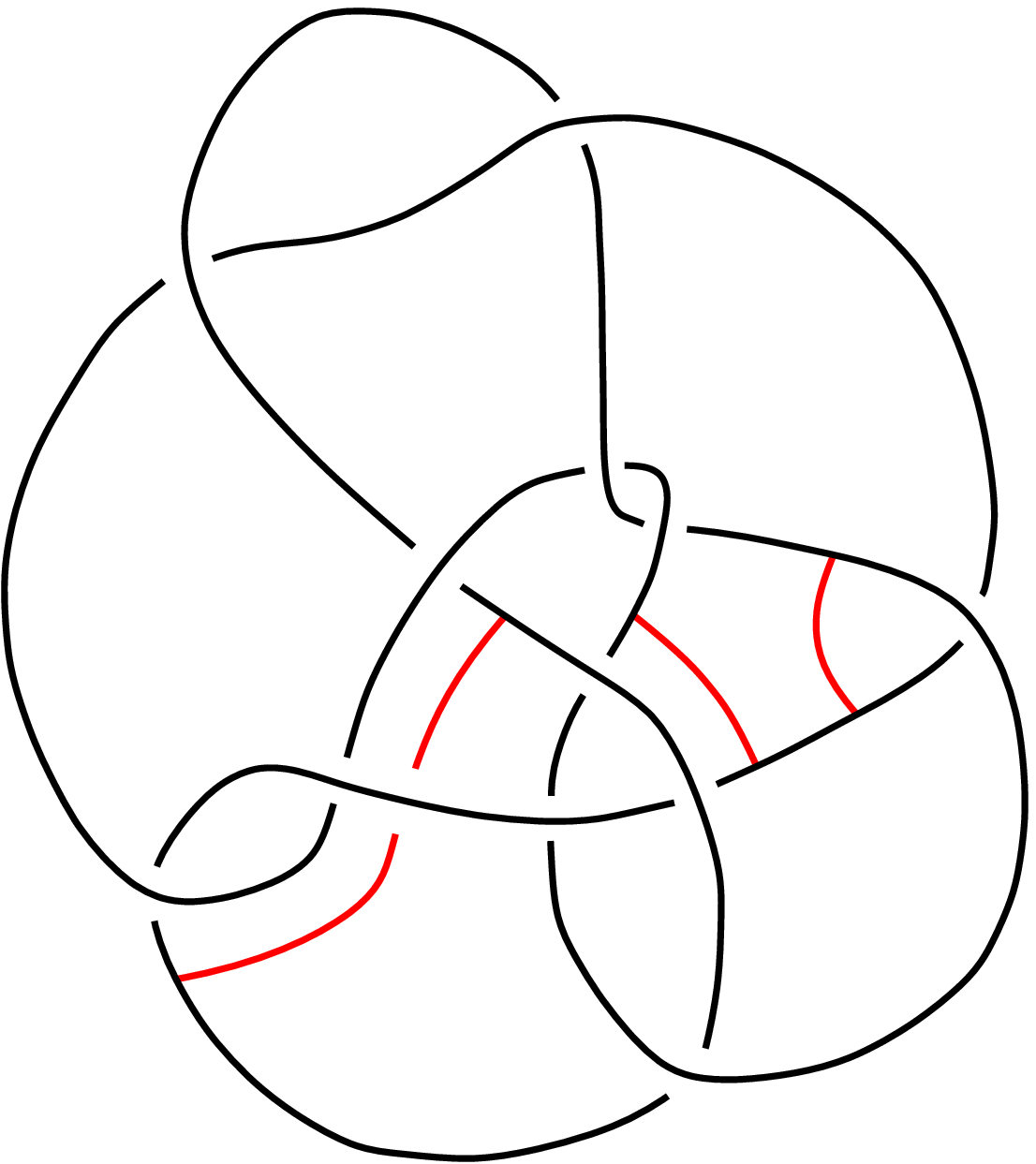}
        \caption{$12n_{769}$}
    \end{subfigure}
    \caption{Genus one slice surfaces}
    \label{fig:genusone}
\end{figure*}

In Figure~\ref{fig:genusone} we present explicit genus one slice surfaces for the knots we discuss in Section~\ref{sec:furtherapplications}.  These surfaces can be constructed from the given diagrams by performing band surgeries along the red arcs in each diagram, where all bands are given the blackboard framing.  Performing these band surgeries yield unlinks, which can be capped off with disks to yield the required genus one slice surfaces.  The only exception is the knot $11a_{37}$, which is transformed into the slice knot $6_1$ by the indicated band surgeries.

Finally, in Tables~\ref{tab:quasipositivepredictions} and \ref{tab:genuspredictions} we present the predictions made by our models for the remaining 11 and 12--crossing knots with unknown quasipositivity and unknown slice genus.  In each case we retrain our original networks with training data that includes the newly discovered quasipositive knots in Table~\ref{tab:quasipositive}, and knots with newly determined slice genus from \cite{Mccoy}.  For each knot type $K$ in Tables~\ref{tab:quasipositivepredictions} and \ref{tab:genuspredictions} predictions were made on all 32 data vectors $v$ representing $K$, with the mean and standard deviation of these 32 predictions presented.  In Table~\ref{tab:quasipositivepredictions} values near 1 indicate knots which are predicted to be quasipositive, while values near 0 correspond to knots which are predicted to not be quasipositive.  Indeed, the values in Table~\ref{tab:quasipositivepredictions} may be interpreted as the estimated probability that the knot $K$ is quasipositive.  

\newpage
\renewcommand{\arraystretch}{1.25}

\begin{center}
\begin{longtable}[c]{  c  c  c  c }
\caption{Knots with quasipositive braid representatives}
\label{tab:quasipositive}\\ \toprule
 {\sc Knot type} & {\sc Quasipositive braid} & {\sc Band centers} & $\tau(K) = g_4(K)$\\ \toprule

 \endfirsthead 
 
\caption*{{\sc Table \ref{tab:quasipositive} Continued:} Knots with quasipositive braid representatives}
\label{tab:quasipositive}\\ \toprule
 {\sc Knot type} & {\sc Quasipositive braid} & {\sc Band centers} & $\tau(K) = g_4(K)$\\ \toprule
 \endhead
 \bottomrule
 \endfoot
 \bottomrule
 \endlastfoot
$11n_{35}$ & $A A b A b A C B B C C$ & 1, 2, 4, 6, 7, 10, 11& 2 \\ \hline
$11n_{40}$ & $A b A b b A C B B B C$ & 1, 3, 6, 7, 11& 1 \\ \hline
$11n_{43}$ & $A A b A b A C B B B C$ & 1, 2, 4, 6, 7, 8, 11& 2 \\ \hline
$11n_{54}$ & $A A B a a B A C b C C$ & 3, 7, 8, 10, 11& 1 \\ \hline
$11n_{59}$ & $A A A b A b C B B B C$ & 1, 2, 3, 5, 7, 8, 11& 2 \\ \hline
$11n_{63}$ & $A A B a B A C b C D c D$ & 1, 2, 5, 7, 10, 12& 1 \\ \hline
$11n_{72}$ & $A A b b A C B B B C C$ & 1, 2, 5, 6, 7, 10, 11& 2 \\ \hline
$11n_{95}$ & $A A A B a B C B a B C$ & 1, 4, 6, 7, 8, 10, 11& 2 \\ \hline
$11n_{105}$ & $A A b A b C B B B C C$ & 1, 2, 4, 6, 7, 10, 11& 2 \\ \hline
$11n_{118}$ & $A A A B A C b A b C B$ & 1, 2, 3, 5, 6, 8, 10& 2 \\ \hline
$11n_{139}$ & $A A B a C b C D c B c D$ & 1, 3, 8, 12& 0 \\ \hline
$11n_{144}$ & $A A b A A C B a B C C$ & 1, 2, 4, 6, 7, 10, 11& 2 \\ \hline
$11n_{162}$ & $A A B a B A C b a C B D c D$ & 1, 5, 7, 11, 12, 14& 1 \\ \hline
$11n_{174}$ & $A A b A C b A C B B C$ & 1, 2, 4, 5, 7, 8, 11& 2 \\ \hline
$11n_{185}$ & $A b A C b A C B B B C$ & 1, 3, 4, 6, 7, 8, 11& 2 \\ \hline
$12n_{79}$ & $A b A B C d C d C D B D B c$ & 1, 3, 7, 9, 11, 13& 1 \\ \hline
$12n_{81}$ & $a B B A A A B c B B a C C$ & 2, 3, 4, 7, 9, 10, 12& 2 \\ \hline
$12n_{117}$ & $A b A B B B c c B B C C C$ & 1, 3, 4, 5, 9, 10, 11& 2 \\ \hline
$12n_{123}$ & $A b A b C C B B c D C C D D$ & 1, 3, 5, 6, 10, 11, 13, 14& 2 \\ \hline
$12n_{128}$ & $A b A b C C D B D B C C D c$ & 1, 3, 5, 6, 7, 9, 11, 13& 2 \\ \hline
$12n_{155}$ & $A B B a B C b b C C B B A$ & 1, 2, 3, 5, 6, 9, 10& 2 \\ \hline
$12n_{157}$ & $A b C D A B B c D B C b$ & 1, 3, 4, 5, 9, 10& 1 \\ \hline
$12n_{176}$ & $A B B A A b C b C D c B c D$ & 1, 2, 4, 5, 10, 14& 1 \\ \hline
$12n_{183}$ & $A b C D A B D B c D B C$ & 1, 3, 4, 5, 6, 7, 10, 11& 2 \\ \hline
$12n_{194}$ & $A A A b C A C A B c B a B$ & 1, 2, 3, 5, 6, 9, 11& 2 \\ \hline
$12n_{209}$ & $a B C C B A A A B c B B a$ & 2, 3, 5, 6, 9, 11, 12& 2 \\ \hline
$12n_{213}$ & $a B C C B B A A A B c B a$ & 2, 3, 5, 6, 7, 10, 12& 2 \\ \hline
$12n_{222}$ & $A B B C A C D A B c B D c b$ & 1, 2, 5, 7, 8, 9, 11, 12& 2 \\ \hline
$12n_{237}$ & $A b b C C B B B B B A b C$ & 1, 4, 5, 6, 7, 11, 13& 2 \\ \hline
$12n_{240}$ & $A A A b b C C B B B A b C$ & 1, 2, 3, 6, 7, 11, 13& 2 \\ \hline
$12n_{249}$ & $A b b A A B B B C D b c D C$ & 1, 4, 5, 9, 10, 13& 1 \\ \hline
$12n_{254}$ & $A B B C A A b b C C B A b$ & 1, 4, 5, 6, 9, 10, 12& 2 \\ 
$12n_{303}$ & $A b A B c c B B C C C C C$ & 1, 3, 7, 8, 9, 10, 11& 2 \\ \hline 
$12n_{306}$ & $A b A B c c B B C C C D c D$ & 1, 3, 7, 8, 12, 14& 1 \\ \hline
$12n_{316}$ & $A b A B c B B c B B C C C$ & 1, 3, 6, 7, 9, 10, 11& 2 \\ \hline
$12n_{321}$ & $a B C A C B a B C A b A B$ & 2, 3, 4, 5, 8, 9, 13& 2 \\ \hline
$12n_{372}$ & $A A A B B C a C a B a B c$ & 4, 5, 6, 10, 12& 1 \\ \hline
$12n_{373}$ & $A B B C a C B B A A b c A$ & 1, 2, 3, 4, 7, 9, 10& 2 \\ \hline
$12n_{375}$ & $A b A B C b C B B B C b C$ & 1, 3, 4, 5, 7, 11, 13& 2 \\ \hline
$12n_{381}$ & $A b C D A B D c B c D C$ & 1, 4, 5, 6, 7, 11& 1 \\ \hline
$12n_{383}$ & $A b C D A A B c B D C d$ & 1, 3, 5, 6, 9, 10& 1 \\ \hline
$12n_{407}$ & $A b A B C b C B c B C C C$ & 1, 3, 4, 7, 11, 12, 13& 2 \\ \hline
$12n_{441}$ & $A A b c A b A B B C C B B$ & 1, 2, 5, 7, 8, 10, 12& 2 \\ \hline
$12n_{487}$ & $a b C C B B A A B C b a B$ & 3, 4, 5, 10, 13& 1 \\ \hline
$12n_{496}$ & $A b C D A B c B D B C C$ & 1, 3, 4, 5, 8, 9, 10, 11& 2 \\ \hline
$12n_{513}$ & $a B A A A B C C C B a B c$ & 2, 3, 6, 7, 8, 10, 12& 2 \\ \hline
$12n_{577}$ & $A b A B C b C B B B c B c$ & 1, 3, 8, 9, 10& 1 \\ \hline
$12n_{582}$ & $A b A B C D b D C B c d B c$ & 1, 3, 6, 10& 0 \\ \hline
$12n_{589}$ & $A b A B C C C B a B c B B$ & 1, 4, 5, 6, 8, 10, 12& 2 \\ \hline
$12n_{671}$ & $A b b A A A B C B B A C C$ & 1, 4, 5, 6, 7, 8, 11, 12, 13& 3 \\ \hline
$12n_{677}$ & $a B C A C A A A B c B a B$ & 2, 3, 4, 6, 9, 11, 13& 2 \\ \hline
$12n_{682}$ & $A A b b A A B C B B A C C$ & 1, 2, 5, 6, 7, 8, 11, 12, 13& 3 \\ \hline
$12n_{719}$ & $A b A B c B c B c B C C C$ & 1, 3, 6, 8, 10& 1 \\ \hline
$12n_{724}$ & $A b C A C A A A B c B a B$ & 1, 3, 4, 6, 7, 9, 11& 2 \\ \hline
$12n_{726}$ & $a B c B C b C A C D c A B D$ & 2, 5, 8, 10, 13, 14& 1 \\ \hline
$12n_{729}$ & $a B C B C A A B c B a B C$ & 2, 3, 4, 5, 8, 10, 12& 2 \\ \hline
$12n_{734}$ & $A B C C B a B B B c a B A$ & 2, 3, 5, 7, 8, 9, 12& 2 \\ \hline
$12n_{735}$ & $A b A B C C D D C B c B C d b c$ & 1, 3, 5, 7, 9, 12& 1 \\ \hline
$12n_{738}$ & $a B B A A B B C C B a B c$ & 2, 3, 6, 7, 8, 10, 12& 2 \\ \hline
$12n_{749}$ & $A B B A A B B A A b b b$ & 1, 2, 4, 5, 8, 9& 2 \\ \hline
$12n_{753}$ & $a B A A B c B a B B C C B$ & 2, 5, 7, 9, 10, 11, 13& 2 \\ \hline
$12n_{770}$ & $a B C C B B A A B c B a B$ & 2, 3, 5, 6, 9, 11, 13& 2 \\ \hline
$12n_{796}$ & $A B B B C C B a B c a B A$ & 2, 3, 4, 5, 7, 9, 12& 2 \\ \hline
$12n_{797}$ & $a B c B C b D c A D C A B c B C$ & 4, 7, 9, 10, 13, 15& 1 \\ 
$12n_{801}$ & $A B a B C b A C B c B C C$ & 1, 2, 5, 8, 9, 11, 12& 2 \\ \hline
$12n_{807}$ & $A b A B C B a C B c B C C$ & 1, 4, 5, 6, 8, 11, 12& 2 \\ \hline
$12n_{811}$ & $a B A b C b A C B c B C C$ & 2, 7, 8, 12, 13& 1 \\ \hline
$12n_{830}$ & $A b b A A B B B A A B B$ & 1, 4, 5, 6, 7, 8, 9, 10& 3 \\ \hline
$12n_{836}$ & $A A B B c a B a B B C C B$ & 3, 4, 7, 9, 10, 11, 13& 2 \\ \hline
$12n_{838}$ & $A B c D a B C A b d C D$ & 1, 4, 6, 11& 0 \\ \hline
$12n_{849}$ & $a B C C B A A A B c B a B$ & 2, 3, 5, 6, 9, 11, 13& 2 \\ \hline
$12n_{863}$ & $A B C C B a B a B B c B A$ & 2, 3, 5, 7, 9, 10, 12& 2 \\ 
\end{longtable}
\end{center}

\begin{table}[h]
\begin{tabular}{  c  c c  c  }
\toprule
{\sc Knot type} & {\sc Quasinegative braid} & {\sc Band centers} & $\tau (K) = -g_4(K)$ \\ 
\toprule
$11n_{1}$ & $a a a b A c B c b d C d$ & 1, 2, 4, 8, 10, 12 & $-1$ \\ \hline
$11n_{10}$ & $a a a b A b b a c B c$ & 1, 2, 3, 6, 7, 9, 11 & $-2$ \\ \hline
$11n_{14}$ & $a a a a b A b a c B c$ & 1, 2, 3, 4, 7, 9, 11 & $-2$ \\ \hline
$11n_{75}$ & $a a B B B a c b b b c$ & 1, 2, 6, 7, 11 & $-1$ \\ \hline
$11n_{84}$ & $a a B a B c b A b b c$ & 1, 2, 6, 9, 11 & $-1$ \\ \hline
$11n_{87}$ & $a a b A b a c B a c B$ & 1, 2, 6, 7, 10 & $-1$ \\ \hline
$11n_{89}$ & $a a a B a c b A b b c$ & 1, 2, 3, 6, 7, 9, 11 & $-2$ \\ \hline
$11n_{108}$ & $a a B a a c b A b b c$ & 1, 2, 4, 6, 7, 9, 11 & $-2$ \\ \hline
$11n_{109}$ & $a a a B a a c b A b c$ & 1, 2, 3, 5, 7, 8, 11 & $-2$ \\ \hline
$11n_{122}$ & $a a a B a B c b A b c$ & 1, 2, 3, 7, 11 & $-1$ \\ \hline
$11n_{134}$ & $a a b A b b c B a B c$ & 1, 2, 5, 7, 11 & $-1$ \\ \hline
$11n_{176}$ & $a a B a B a c b A b c$ & 1, 2, 6, 7, 11 & $-1$ \\ \bottomrule
\end{tabular}
\caption{Knots with quasinegativebraid representatives}
\label{tab:quasinegative}
\end{table}

\FloatBarrier

\begin{center}
\begin{longtable}[c]{  c  c  c  }
\caption{Knots with mean quasipositivity predictions}
\label{tab:quasipositivepredictions}\\ \toprule

 {\sc Knot type} & {\sc Mean quasipositive probability} & {\sc Standard deviation} \\ 
\toprule
 \endfirsthead 
\caption*{{\sc Table~\ref{tab:quasipositivepredictions} Continued:} Knots with mean quasipositivity predictions}\\
\toprule
 {\sc Knot type} & {\sc Mean quasipositive probability} & {\sc Standard deviation} \\ 
\toprule
 \endhead
 \bottomrule
 \endfoot
 \bottomrule
 \endlastfoot
$11n_{17}$ & $0.15625$ & $0.368902033$ \\ \hline
$11n_{22}$ & $0.561580203$ & $0.503215897$ \\ \hline
$11n_{37}$ & $1.38\times 10^{-21}$ & $7.80\times 10^{-21}$ \\ \hline
$11n_{46}$ & $0.28250862$ & $0.456058364$ \\ \hline
$11n_{50}$ & $0.1871539$ & $0.377428619$ \\ \hline
$11n_{71}$ & $0.468751151$ & $0.50700625$ \\ \hline
$11n_{91}$ & $0.25$ & $0.439941345$ \\ \hline
$11n_{99}$ & $0.906253413$ & $0.296133755$ \\ \hline
$11n_{113}$ & $0.03125$ & $0.176776695$ \\ \hline
$11n_{127}$ & $0.970651054$ & $0.166022555$ \\ \hline
$11n_{132}$ & $0.366879124$ & $0.482029503$ \\ \hline
$11n_{146}$ & $0.243476685$ & $0.426402697$ \\ \hline
$11n_{159}$ & $0.767814315$ & $0.416174674$ \\ \hline
$11n_{172}$ & $0.294072376$ & $0.440207634$ \\ \hline
$11n_{178}$ & $0.406002608$ & $0.492858534$ \\ \hline
$11n_{184}$ & $0.344561608$ & $0.48197181$ \\ \hline
$12n_{5}$ & $1.69\times 10^{-28}$ & $9.55\times 10^{-28}$ \\ \hline
$12n_{58}$ & $0$ & $0$ \\ \hline
$12n_{72}$ & $0.375$ & $0.491869377$ \\ \hline
$12n_{76}$ & $0.995966416$ & $0.022817398$ \\ \hline
$12n_{80}$ & $0.062509919$ & $0.245932093$ \\ \hline
$12n_{116}$ & $9.38\times 10^{-17}$ & $5.30\times 10^{-16}$ \\ \hline
$12n_{120}$ & $0.674233288$ & $0.465984477$ \\ \hline
$12n_{121}$ & $0.587734634$ & $0.495018156$ \\ \hline
$12n_{140}$ & $0.03125$ & $0.176776695$ \\ \hline
$12n_{145}$ & $3.28\times 10^{-11}$ & $1.86\times 10^{-10}$ \\ \hline
$12n_{146}$ & $0$ & $0$ \\ \hline
$12n_{148}$ & $0.906249916$ & $0.296144554$ \\ \hline
$12n_{149}$ & $9.30\times 10^{-30}$ & $5.26\times 10^{-29}$ \\ \hline
$12n_{159}$ & $0.162582742$ & $0.36696125$ \\ \hline
$12n_{168}$ & $1$ & $0$ \\ \hline
$12n_{171}$ & $0.90625$ & $0.296144581$ \\ 
$12n_{193}$ & $0.220194975$ & $0.413602013$ \\ \hline
$12n_{199}$ & $0.09375$ & $0.296144581$ \\ \hline
$12n_{200}$ & $0$ & $0$ \\ \hline
$12n_{208}$ & $0.113261379$ & $0.307961737$ \\ \hline
$12n_{212}$ & $0.156250003$ & $0.368902031$ \\ \hline
$12n_{236}$ & $0.03125$ & $0.176776695$ \\ \hline
$12n_{239}$ & $0.031209243$ & $0.176546141$ \\ \hline
$12n_{247}$ & $0.875$ & $0.336010753$ \\ \hline
$12n_{253}$ & $0.156339736$ & $0.368863146$ \\ \hline
$12n_{260}$ & $0$ & $0$ \\ \hline
$12n_{270}$ & $0.03125001$ & $0.176776694$ \\ \hline
$12n_{290}$ & $1$ & $0$ \\ \hline
$12n_{293}$ & $1$ & $0$ \\ \hline
$12n_{312}$ & $0$ & $0$ \\ \hline
$12n_{318}$ & $6.60\times 10^{-9}$ & $3.74\times 10^{-8}$ \\ \hline
$12n_{332}$ & $1.27\times 10^{-10}$ & $7.17\times 10^{-10}$ \\ \hline
$12n_{347}$ & $0.093749904$ & $0.296144275$ \\ \hline
$12n_{366}$ & $0.90625$ & $0.296144581$ \\ \hline
$12n_{379}$ & $2.13\times 10^{-9}$ & $7.33\times 10^{-9}$ \\ \hline
$12n_{393}$ & $2.75\times 10^{-36}$ & $1.56\times 10^{-35}$ \\ \hline
$12n_{397}$ & $0$ & $0$ \\ \hline
$12n_{404}$ & $3.14\times 10^{-22}$ & $1.78\times 10^{-21}$ \\ \hline
$12n_{409}$ & $0.104981924$ & $0.299270629$ \\ \hline
$12n_{414}$ & $4.31\times 10^{-40}$ & $2.44\times 10^{-39}$ \\ \hline
$12n_{429}$ & $0$ & $0$ \\ \hline
$12n_{432}$ & $0.307731149$ & $0.463982403$ \\ \hline
$12n_{451}$ & $0.2188284$ & $0.41997143$ \\ \hline
$12n_{454}$ & $0.533521886$ & $0.504434486$ \\ \hline
$12n_{469}$ & $0.185870042$ & $0.393203451$ \\ \hline
$12n_{510}$ & $0.03125$ & $0.176776695$ \\ \hline
$12n_{512}$ & $0.0625$ & $0.245934688$ \\ \hline
$12n_{514}$ & $0.031250241$ & $0.176776651$ \\ 
$12n_{520}$ & $0.03125$ & $0.176776695$ \\ \hline
$12n_{522}$ & $0.374219347$ & $0.474985581$ \\ \hline
$12n_{523}$ & $0$ & $0$ \\ \hline
$12n_{528}$ & $0.906251277$ & $0.296140546$ \\ \hline
$12n_{543}$ & $0.820412473$ & $0.381819891$ \\ \hline
$12n_{549}$ & $0.590909394$ & $0.496498825$ \\ \hline
$12n_{564}$ & $0$ & $0$ \\ \hline
$12n_{572}$ & $0.013379612$ & $0.075686466$ \\ \hline
$12n_{606}$ & $0.125000239$ & $0.336010632$ \\ \hline
$12n_{621}$ & $0$ & $0$ \\ \hline
$12n_{626}$ & $0$ & $0$ \\ \hline
$12n_{642}$ & $0.172397691$ & $0.363832949$ \\ \hline
$12n_{660}$ & $0.999999981$ & $1.05\times 10^{-7}$ \\ \hline
$12n_{667}$ & $0.426719756$ & $0.495188198$ \\ \hline
$12n_{685}$ & $1.20\times 10^{-30}$ & $6.59\times 10^{-30}$ \\ \hline
$12n_{698}$ & $0.399001675$ & $0.491717623$ \\ \hline
$12n_{699}$ & $1.06\times 10^{-20}$ & $6.00\times 10^{-20}$ \\ \hline
$12n_{700}$ & $0.216394423$ & $0.415638984$ \\ \hline
$12n_{701}$ & $0.031228665$ & $0.176656008$ \\ \hline
$12n_{717}$ & $0.284256679$ & $0.443798571$ \\ \hline
$12n_{730}$ & $0.692569353$ & $0.46410254$ \\ \hline
$12n_{742}$ & $0.125$ & $0.336010753$ \\ \hline
$12n_{768}$ & $1.01\times 10^{-7}$ & $5.74\times 10^{-7}$ \\ \hline
$12n_{769}$ & $0.187499929$ & $0.39655762$ \\ \hline
$12n_{771}$ & $0.21986673$ & $0.419459829$ \\ \hline
$12n_{814}$ & $0.593775695$ & $0.498959377$ \\ \hline
$12n_{823}$ & $0.250065521$ & $0.438135766$ \\ \hline
$12n_{861}$ & $0.242941631$ & $0.429199989$ \\ \hline
$12n_{862}$ & $0.364586726$ & $0.481107995$ \\ \hline
$12n_{867}$ & $0.312510141$ & $0.47092213$ \\ \hline
$12n_{871}$ & $1.22\times 10^{-27}$ & $6.89\times 10^{-27}$ \\ 
\end{longtable}
\end{center}

\begin{table}[h]
\begin{tabular}{  c  c  c  }
\toprule
{\sc Knot type} & {\sc Mean slice genus prediction} & {\sc Standard deviation} \\ 
\toprule
$11n_{34}$ & $0.496090381$ & $0.303309781$ \\ \hline
$11n_{80}$ & $1.035585403$ & $0.039595998$ \\ \hline
$12a_{153}$ & $1.069947992$ & $0.081041807$ \\ \hline
$12a_{187}$ & $1.009274485$ & $0.020159284$ \\ \hline
$12a_{230}$ & $1.019275038$ & $0.029479206$ \\ \hline
$12a_{317}$ & $1.030372641$ & $0.06259502$ \\ \hline
$12a_{450}$ & $1.009170415$ & $0.013387401$ \\ \hline
$12a_{570}$ & $1.009285686$ & $0.012420726$ \\ \hline
$12a_{624}$ & $1.027900245$ & $0.042808674$ \\ \hline
$12a_{636}$ & $1.005604245$ & $3.13\times 10^{-5}$ \\ \hline
$12a_{787}$ & $1.005587707$ & $3.58\times 10^{-5}$ \\ \hline
$12a_{905}$ & $1.005527869$ & $3.22\times 10^{-5}$ \\ \hline
$12a_{1189}$ & $1.07924728$ & $0.092337155$ \\ \hline
$12a_{1208}$ & $1.023944944$ & $0.028467024$ \\ \hline
$12n_{52}$ & $1.024184935$ & $0.040119993$ \\ \hline
$12n_{63}$ & $1.010492662$ & $0.066061094$ \\ \hline
$12n_{225}$ & $1.001823799$ & $0.028759343$ \\ \hline
$12n_{239}$ & $1.016089445$ & $0.025106408$ \\ \hline
$12n_{269}$ & $1.00551357$ & $8.59\times 10^{-6}$ \\ \hline
$12n_{505}$ & $1.037368491$ & $0.074247022$ \\ \hline
$12n_{512}$ & $1.027464205$ & $0.049029925$ \\ \hline
$12n_{542}$ & $1.038927461$ & $0.057118994$ \\ \hline
$12n_{555}$ & $1.022189596$ & $0.022822694$ \\ \hline
$12n_{558}$ & $1.013406007$ & $0.02690478$ \\ \hline
$12n_{598}$ & $1.0063358$ & $0.004468055$ \\ \hline
$12n_{602}$ & $1.010329522$ & $0.012901433$ \\ \hline
$12n_{665}$ & $1.006596039$ & $0.029156534$ \\ \hline
$12n_{756}$ & $1.005524328$ & $7.54\times 10^{-5}$ \\ \hline
$12n_{886}$ & $1.006042294$ & $0.002597384$ \\ \bottomrule
\end{tabular}
\caption{Knots with mean slice genus predictions}
\label{tab:genuspredictions}
\end{table}

\clearpage

\bibliographystyle{plain}
\bibliography{bibliography}
\end{document}